\documentclass[12pt,a4paper]{article}
\usepackage{amsfonts,amssymb,amsmath,indentfirst,amsthm}

\input{tcilatex}

\begin{document}

\date{}
\title{Some New Inequalities of Dirichlet Eigenvalues for Laplace Operator with any Order}
\author{Na Huang \thanks{%
Na Huang: Department of Applied Mathematics; Key Laboratory of Space
Applied Physics and Chemistry, Ministry of Education, Northwestern
Polytechnical University, Xi'an, Shaanxi, 710129, China. e-mail:
huangna7@126.com},
Pengcheng Niu \thanks{%
Pengcheng Niu: Corresponding author, Department of Applied
Mathematics; Key Laboratory of Space Applied Physics and Chemistry,
Ministry of Education, Northwestern Polytechnical University, Xi'an,
Shaanxi, 710129, China. e-mail: pengchengniu@nwpu.edu.cn}
}
\date{}
\maketitle

\newtheorem{thm}{{\indent}Theorem} \newtheorem{cor}[thm]{{\indent}Corollary}
\newtheorem{lem}[thm]{{\indent}Lemma} \newtheorem{prop}[thm]{%
{\indent}Proposition} \theoremstyle{definition} \newtheorem{defn}[thm]{%
{\indent}Definition}

\newtheorem{theo}{\hspace*{2.0em}Theorem\hspace*{0.1em}}[section]
\newenvironment{keywords}{\par\textbf{keywords:}\mbox{ }}{ } \newenvironment{%
coj}{\par\textbf{Conjecture:}\mbox{ }}{ } \numberwithin{equation}{section}

\textbf{Abstract.} In this paper, we establish several inequalities
of Dirichlet eigenvalues for Laplace operator $\Delta $ with any
order on \emph{n}-dimensional Euclidean space. These inequalities
are more general than known Yang's inequalities and contain new
consequences. To obtain them, we borrow the approach of Illias and
Makhoul, and use a generalized Chebyshev's inequality.


\section{Introduction}

\label{introduction}

The following Dirichlet problem
\begin{equation}
\left\{
\begin{array}{l}
{\left( { - \Delta } \right)^l}u = \lambda u,{\kern 1pt} {\kern 1pt} {\kern 1pt} {\kern 1pt} {\kern 1pt} {\kern 1pt} {\kern 1pt} {\kern 1pt} {\kern 1pt} {\kern 1pt} {\kern 1pt} {\kern 1pt} {\kern 1pt} {\kern 1pt} {\kern 1pt} {\kern 1pt} {\kern 1pt} {\kern 1pt} {\kern 1pt} {\kern 1pt} {\kern 1pt} {\kern 1pt} {\kern 1pt} {\kern 1pt} {\kern 1pt} {\kern 1pt} {\kern 1pt} {\kern 1pt} {\kern 1pt} {\kern 1pt} {\kern 1pt} {\kern 1pt} {\kern 1pt} {\kern 1pt} {\kern 1pt} {\kern 1pt} {\kern 1pt} {\kern 1pt} {\kern 1pt} {\kern 1pt} {\kern 1pt} {\kern 1pt} {\kern 1pt} {\kern 1pt} {\kern 1pt} {\kern 1pt} {\kern 1pt} {\kern 1pt} {\kern 1pt} {\kern 1pt} {\kern 1pt} {\kern 1pt} {\kern 1pt} {\kern 1pt} {\kern 1pt} {\kern 1pt} {\kern 1pt} {\kern 1pt} {\kern 1pt} {\kern 1pt} {\kern 1pt} {\kern 1pt} {\kern 1pt} {\kern 1pt} {\kern 1pt} {\kern 1pt} {\kern 1pt} \,\,\,{\rm{in}} \,\,\,\Omega ,\\
u = \frac{{\partial u}}{{\partial \nu }} =  \cdot  \cdot  \cdot  =
\frac{{{\partial ^{l - 1}}u}}{{\partial {\nu ^{l - 1}}}} = 0,{\kern
1pt} {\kern 1pt} {\kern 1pt} {\kern 1pt} {\kern 1pt} {\kern 1pt}
{\kern 1pt} {\kern 1pt} {\kern 1pt} {\kern 1pt} {\kern 1pt} {\kern
1pt} {\kern 1pt} {\rm{on}} \,\,\partial \Omega
\end{array} \right.
\end{equation}
has been extensively considered, where $\Delta $ is the Laplacian:
$\Delta  = \sum\limits_{j = 1}^n {\frac{{{\partial ^2}}}{{\partial
x_j^2}}} $, $l$ is any positive integer, $\Omega$ is a bounded
domain in the Euclidean space ${R^n}$, $\nu $ is the outward unit
normal on $\partial \Omega $.

When $l = 1$, Payne, P\'{o}lya and Weinberger in [9] showed the
following inequality (the PPW inequality)
\[{\lambda _{k + 1}} - {\lambda _k} \le \frac{4}{{nk}}\sum\limits_{r = 1}^k {{\lambda
_r}}.
\]
The inequality (the HP inequality)
\[\sum\limits_{r = 1}^k {\frac{{{\lambda _r}}}{{{\lambda _{k + 1}} - {\lambda _r}}}}  \ge
\frac{{nk}}{4}\] was due to Hile and Protter in [4]. Yang in [10]
proved some important eigenvalue estimates, which are Yang's first
inequality
\[\sum\limits_{r = 1}^k {{{\left( {{\lambda _{k + 1}} - {\lambda _r}} \right)}^2}}  \le \frac{4}{n}\sum\limits_{r = 1}^k {\left( {{\lambda _{k + 1}} - {\lambda _r}} \right)} {\lambda
_r}\] and Yang's second inequality
\[{\lambda _{k + 1}} \le \left( {1 + \frac{4}{n}} \right)\frac{1}{k}\sum\limits_{r = 1}^k {{\lambda _r}}
.\]

When $l = 2$, the estimate
\[{\lambda _{k + 1}} - {\lambda _k} \le \frac{{8(n + 2)}}{{{n^2}k}}\sum\limits_{r = 1}^k {{\lambda _r}} \]
was derived by Payne, P¨®lya and Weinberger in [9]. Chen and Qian in
[1] and Hook in [5] proved respectively
\[\frac{{{n^2}{k^2}}}{{8(n + 2)}} \le \sum\limits_{i = 1}^k {\frac{{\lambda _i^{\frac{1}{2}}}}{{{\lambda _{k + 1}} - {\lambda _i}}}} \sum\limits_{i = 1}^k {\lambda _i^{\frac{1}{2}}} .\]
The following inequality
\[\sum\limits_{i = 1}^k {({\lambda _{k + 1}} - {\lambda _i})}  \le {\left( {\frac{{8(n + 2)}}{{{n^2}}}} \right)^{\frac{1}{2}}}\sum\limits_{i = 1}^k {{{\left( {{\lambda _i}({\lambda _{k + 1}} - {\lambda _i})} \right)}^{\frac{1}{2}}}} \]
was gotten by Cheng and Yang [3].

For any positive integer $l$, Chen and Qian [1] and Hook [5]
independently obtained
\[\frac{{{n^2}{k^2}}}{{4l\left( {n + 2l - 2} \right)}} \le \sum\limits_{i = 1}^k {\frac{{\lambda _i^{\frac{1}{l}}}}{{{\lambda _{k + 1}} - {\lambda _i}}}} \sum\limits_{i = 1}^k {\lambda _i^{\frac{{l - 1}}{l}}} ,k = 1,2, \cdot  \cdot  \cdot .\]
The inequality
\[\sum\limits_{i = 1}^k {({\lambda _{k + 1}} - {\lambda _i}} {)^2} \le \frac{{4l(n + 2l - 2)}}{{{n^2}}}\sum\limits_{i = 1}^k {({\lambda _{k + 1}} - {\lambda _i})} {\lambda _i}\]
was concluded by Cheng, Ichikawa and Mametsuka in [2]. Ilias and
Makhoul in [7] exhibited a new abstract formula relating eigenvalues
of a self-adjoint operator and deduced Yang type inequality for
Dirichlet eigenvalues of sub-Laplacian with any order on the
Heisenberg group, in the light of Chebyshev's inequality.

In this paper, we will give several new estimates of Dirichlet
eigenvalues to (1.1) by combining the approach of Ilias and Makhoul
in [8] and using a generalized Chebyshev's inequality in [6]. For
convenience, we denote $L =  - \Delta $ and assume always that
${\lambda _{i + 1}} > {\lambda _i},i = 1,2, \cdot  \cdot  \cdot ,$
in the sequel. The main results of this paper are the following
Theorem 1.1 and its corollaries.

\textbf{Theorem 1.1} Let $\left\{ {{\lambda _i}} \right\}$ be the
eigenvalues of (1.1), then
\begin{align}
&\sum\limits_{i = 1}^k {({\lambda _{k + 1}} - {\lambda _i}}
{)^\alpha } \notag
\\
&\le  \frac{{2\sqrt {l(n + 2l - 2)} }}{n}{\left[ {\sum\limits_{i =
1}^k {{{({\lambda _{k + 1}} - {\lambda _i})}^\beta }\lambda
_i^{\frac{{l - 1}}{l}}} } \right]^{\frac{1}{2}}}{\left[
{\sum\limits_{i = 1}^k {{{({\lambda _{k + 1}} - {\lambda
_i})}^{2\alpha  - \beta  - 1}}\lambda _i^{\frac{1}{l}}} }
\right]^{\frac{1}{2}}},
\end{align}
where $\alpha  \in R$ and $\beta  \ge 0$ such that ${\alpha ^2} \le
2\beta $.

\textbf{Remark 1.1} Inequality (1.2) is the generalization of Yang's
inequality. Some consequences are easily deduced from (1.2) and new
inequalities are listed now.

(1) Let $2\alpha  - \beta  - 1 = 0$, then
\[\sum\limits_{i = 1}^k {({\lambda _{k + 1}} - {\lambda _i}} {)^\alpha } \le \frac{{2\sqrt {l(n + 2l - 2)} }}{n}{\left[ {\sum\limits_{i = 1}^k {{{({\lambda _{k + 1}} - {\lambda _i})}^{2\alpha  - 1}}\lambda _i^{\frac{{l - 1}}{l}}} } \right]^{\frac{1}{2}}}{\left[ {\sum\limits_{i = 1}^k {\lambda _i^{\frac{1}{l}}} } \right]^{\frac{1}{2}}},\]
where $\alpha  \in \left[ {2 - \sqrt 2 ,2 + \sqrt 2 } \right]$.

(2) When $\alpha  =\beta= 1$, it follows
\[\sum\limits_{i = 1}^k {({\lambda _{k + 1}} - {\lambda _i}} ) \le \frac{{2\sqrt {l(n + 2l - 2)} }}{n}{\left[ {\sum\limits_{i = 1}^k {({\lambda _{k + 1}} - {\lambda _i})\lambda _i^{\frac{{l - 1}}{l}}} } \right]^{\frac{1}{2}}}{\left[ {\sum\limits_{i = 1}^k {\lambda _i^{\frac{1}{l}}} } \right]^{\frac{1}{2}}}.\]

(3) When $\alpha  = \frac{1}{2}$, and $\beta  \ge \frac{1}{8}$, we
have
\[\sum\limits_{i = 1}^k {({\lambda _{k + 1}} - {\lambda _i}} {)^{\frac{1}{2}}} \le \frac{{2\sqrt {l(n + 2l - 2)} }}{n}{\left[ {\sum\limits_{i = 1}^k {{{({\lambda _{k + 1}} - {\lambda _i})}^\beta }\lambda _i^{\frac{{l - 1}}{l}}} } \right]^{\frac{1}{2}}}{\left[ {\sum\limits_{i = 1}^k {{{({\lambda _{k + 1}} - {\lambda _i})}^{ - \beta }}\lambda _i^{\frac{1}{l}}} } \right]^{\frac{1}{2}}}.\]

(4) When $\alpha  =  - 1$, and $\beta  \ge \frac{1}{2}$, then
\[\sum\limits_{i = 1}^k {({\lambda _{k + 1}} - {\lambda _i}} {)^{ - 1}} \le \frac{{2\sqrt {l(n + 2l - 2)} }}{n}{\left[ {\sum\limits_{i = 1}^k {{{({\lambda _{k + 1}} - {\lambda _i})}^\beta }\lambda _i^{\frac{{l - 1}}{l}}} } \right]^{\frac{1}{2}}}{\left[ {\sum\limits_{i = 1}^k {{{({\lambda _{k + 1}} - {\lambda _i})}^{ - \beta  - 3}}\lambda _i^{\frac{1}{l}}} } \right]^{\frac{1}{2}}}.\]

We note that the forms in (3) and (4) are never seen previously.

\textbf{Corollary 1.1} We have
\begin{equation}
\sum\limits_{i = 1}^k
{({\lambda _{k + 1}} - {\lambda _i}} {)^2} \le \frac{{2\sqrt {l(n +
2l - 2)} }}{n}{\left[ {\sum\limits_{i = 1}^k {{{({\lambda _{k + 1}}
- {\lambda _i})}^2}} \lambda _i^{\frac{{l - 1}}{l}}}
\right]^{\frac{1}{2}}}{\left[ {\sum\limits_{i = 1}^k {({\lambda _{k
+ 1}} - {\lambda _i})} \lambda _i^{\frac{1}{l}}}
\right]^{\frac{1}{2}}}.
\end{equation}

\textbf{Corollary 1.2} It holds
\begin{equation}
{\lambda _{k + 1}} -
{\lambda _k} \le \frac{{4l(n + 2l - 2)}}{{{n^2}{k^2}}}\left(
{\sum\limits_{i = 1}^k {\lambda _i^{\frac{{l - 1}}{l}}} }
\right)\left( {\sum\limits_{i = 1}^k {\lambda _i^{\frac{1}{l}}} }
\right).
\end{equation}

\textbf{Corollary 1.3} We have
\begin{align}
&\sum\limits_{i = 1}^k {({\lambda _{k + 1}} - {\lambda _i}}
{)^\alpha } \notag
\\
&\le \frac{{2\sqrt {l(n + 2l - 2)} }}{n}{\left[ {\sum\limits_{i =
1}^k {{{({\lambda _{k + 1}} - {\lambda _i})}^\beta }} }
\right]^{\frac{1}{2}}}{\left[ {\sum\limits_{i = 1}^k {{{({\lambda
_{k + 1}} - {\lambda _i})}^{2\alpha  - \beta  - 1}}{\lambda _i}} }
\right]^{\frac{1}{2}}},
\end{align}
where $\alpha  \in R,\beta  \ge 0$ and ${\alpha ^2} \le 2\beta $.

\textbf{Corollary 1.4} Yang type first inequality for (1.1) holds:
\begin{equation}
\sum\limits_{i = 1}^k {({\lambda _{k + 1}} - {\lambda _i}}
{)^2} \le \frac{{4l(n + 2l - 2)}}{{{n^2}}}\sum\limits_{i = 1}^k
{({\lambda _{k + 1}} - {\lambda _i})} {\lambda _i}.
\end{equation}

\textbf{Corollary 1.5} We have the Payne-P\'{o}lya-Weinberger Type
inequality \begin{equation}{\lambda _{k + 1}} - {\lambda _k} \le
\frac{{4l(n + 2l - 2)}}{{{n^2}k}}\sum\limits_{i = 1}^k {{\lambda
_i}}.\end{equation}

\textbf{Corollary 1.6} The Yang type second inequality holds:
\begin{equation}
{\lambda _{k + 1}} \le \left( {1 + \frac{{4l(n + 2l - 2)}}{{{n^2}}}}
\right)\frac{1}{k}\left( {\sum\limits_{i = 1}^k {{\lambda _i}} }
\right).\end{equation}

This paper is arranged as follows. Section 2 is devoted to the
description of known results and some elementary inequalities. The
proofs of Theorem 1.1 and Corollaries 1.1-1.6 are given in Section
3.


\section{Preliminaries}

\label{2}

\textbf{Definition 2.1}(see [7]) A couple $(f,g)$ of functions $f$
and $g$ on the interval $(0,\lambda )$ ($\lambda  > 0$) is said to
belong to ${\chi _\lambda }$ provided that

(1) $f$ and $g$ are positive;

(2) $f$ and $g$ satisfy
\[{\left( {\frac{{f(x) - f(y)}}{{x - y}}} \right)^2} + {\left( {\frac{{{{\left( {f(x)} \right)}^2}}}{{g(x)(\lambda  - x)}} + \frac{{{{\left( {f(y)} \right)}^2}}}{{g(y)(\lambda  - y)}}} \right)}\left( {\frac{{g(x) - g(y)}}{{x - y}}} \right) \le 0,\]
for any $x,y \in (0,\lambda )$, $x \ne y$.

\textbf{Lemma 2.1}(see [6]) Let $(f,g) \in {\chi _\lambda }$, then
$g$ must be nonincreasing; if $f(x) = {(\lambda  - x)^\alpha }$,
$g(x) = {(\lambda  - x)^\beta }$, then ${\alpha ^2} \le 2\beta $.

\textbf{Definition 2.2}(see [7]) For any two operators $A$ and $B$,
their commutator $\left[ {A,B} \right]$ is defined by $\left[ {A,B}
\right] = AB - BA$.

\textbf{Lemma 2.2} For $p = 1,2, \cdot  \cdot  \cdot ,n,$ we have
\begin{equation}
{L^l}({x_p}{u_i}) = {x_p}{L^l}{u_i} - 2l{L^{l -
1}}\frac{\partial }{{\partial {x_p}}}{u_i},
\end{equation}
\begin{equation}
\left[ {L^l,{x_p}} \right]{u_i} =  - 2l{L^{l - 1}}\frac{\partial
}{{\partial {x_p}}}{u_i}.
\end{equation}

\textbf{Proof. } When $l = 1$, we have
\[\frac{\partial }{{\partial {x_j}}}\left( {{x_p}{u_i}} \right) = \left( {\frac{\partial }{{\partial {x_j}}}{x_p}} \right){u_i} + {x_p}\left( {\frac{\partial }{{\partial {x_j}}}{u_i}} \right)\]
and
\[\frac{{{\partial ^2}}}{{\partial x_j^2}}({x_p}{u_i}) = 2\left( {\frac{\partial }{{\partial {x_j}}}{x_p}} \right)\left( {\frac{\partial }{{\partial {x_j}}}{u_i}} \right) + {x_p}\left( {\frac{{{\partial ^2}}}{{\partial x_j^2}}{u_i}} \right).\]
Hence
\[L({x_p}{u_i}) = ( - \Delta )({x_p}{u_i}) = {x_p}L{u_i} - 2\frac{\partial }{{\partial {x_p}}}{u_i},\]
and (2.1) is proved.

Assuming (2.1) is true for $l - 1$, direct calculations show
\[{L^{l - 1}}({x_p}{u_i}) = {x_p}{L^{l - 1}}{u_i} - 2(l - 1){L^{l - 2}}\frac{\partial }{{\partial {x_p}}}{u_i}\]
and
\[\begin{array}{l}
{L^l}({x_p}{u_i}) = L({L^{l - 1}}({x_p}{u_i}))\\
 = L({x_p}{L^{l - 1}}{u_i} - 2(l - 1){L^{l - 2}}\frac{\partial }{{\partial {x_p}}}{u_i})\\
 = {x_p}{L^l}{u_i} - 2l{L^{l - 1}}\frac{\partial }{{\partial {x_p}}}{u_i}.
\end{array}\]
So (2.1) is valid for $l$.

Noting
\[\left[ {{L^l},{x_p}} \right]{u_i} = {L^l}({x_p}{u_i}) - {x_p}{L^l}{u_i}  =  - 2l{L^{l - 1}}\frac{\partial }{{\partial {x_p}}}{u_i},\]
it follows (2.2).

\textbf{Lemma 2.3}(see [7]) Let $A$: $D \subset H \to H$ be a
self-adjoint operator defined on a dense domain $D$, which is
semibounded below and has a discrete spectrum ${\lambda _1} \le
{\lambda _2} \le {\lambda _3} \cdot  \cdot  \cdot $. Let $\left\{
{{T_p}:D \to H} \right\}_{p = 1}^n$ be a collection of
skew-symmetric operators and $\left\{ {{B_p}:{T_p}(D) \to H}
\right\}_{p = 1}^n$ a collection of symmetric operators, leaving $D$
invariant. We denote by $\left\{ {{u_i}} \right\}_{i = 1}^n$ a basis
of orthonormal eigenvectors of $A$, ${u_i}$ corresponding to
${\lambda _i}$ and let ${\lambda _{k + 1}} \ge {\lambda _k}$, $k \ge
1$. Then for any $(f,g)$ in ${\chi _{{\lambda _{k + 1}}}}$, it
follows
\begin{align}
&{\left( {\sum\limits_{i = 1}^k
{\sum\limits_{p = 1}^n {f({\lambda _i})\left\langle {\left[
{{T_p},{B_p}} \right]{u_i},{u_i}} \right\rangle } } }
\right)^2}\notag\\
& \le 4\left( {\sum\limits_{i = 1}^k {\sum\limits_{p = 1}^n
{g({\lambda _i})\left\langle {\left[ {A,{B_p}}
\right]{u_i},{B_p}{u_i}} \right\rangle } } } \right)\left(
{\sum\limits_{i = 1}^k {\sum\limits_{p = 1}^n {\frac{{{{(f({\lambda
_i}))}^2}}}{{g({\lambda _i})({\lambda _{k + 1}} - {\lambda
_i})}}{{\left\| {{T_p}{u_i}} \right\|}^2}} } } \right).
\end{align}

\textbf{Lemma 2.4}(see [6]) For $\gamma  \ge 1$, ${s_i} \ge 0,i = 1,
\cdot  \cdot  \cdot ,k$, it follows
\[{\left( {\sum\limits_{i = 1}^k {{s_i}} } \right)^\gamma } \le {k^{\gamma  - 1}}\sum\limits_{i = 1}^k {s_i^\gamma }. \]

\textbf{Lemma 2.5}(Chebyshev's inequality, [8]) If $\left( {{a_k} -
{a_j}} \right)\left( {{b_k} - {b_j}} \right) \le 0$ for any
nonnegative $k,j$, then
\[\sum\limits_{i = 1}^n {{a_i}{b_i}}  \le \frac{1}{n}\left( {\sum\limits_{i = 1}^n {{a_i}} } \right)\left( {\sum\limits_{i = 1}^n {{b_i}} } \right).\]

\textbf{Lemma 2.6}(generalized Chebyshev's inequality, see [6]) If
${A_1} \ge {A_2} \ge \cdot \cdot \cdot  \ge {A_k} \ge 0,$ $0 \le
{B_1} \le {B_2} \le \cdot \cdot \cdot  \le {B_k},$ $0 \le {C_1} \le
{C_2} \le  \cdot \cdot \cdot \le {C_k},$ then it implies that for
${\alpha ^{^2}} \le 2\beta $,
\begin{equation}
\sum\limits_{i = 1}^k {A_i^\beta {B_i}} \sum\limits_{i = 1}^k
{A_i^{2\alpha  - \beta  - 1}{C_i}}  \le \sum\limits_{i = 1}^k
{A_i^\beta } \sum\limits_{i = 1}^k {A_i^{2\alpha  - \beta  -
1}{B_i}{C_i}}.
\end{equation}

By Lemma 2.6, we immediately have

\textbf{Corollary 2.1}(see [2]) If ${A_1} \ge {A_2} \ge  \cdot \cdot
\cdot \ge {A_k} \ge 0,$ $0 \le {B_1} \le {B_2} \le  \cdot \cdot
\cdot \le {B_k},$ $0 \le {C_1} \le {C_2} \le  \cdot  \cdot \cdot \le
{C_k},$ then we have that for ${\alpha ^{^2}} \le 2\beta $,
\[\sum\limits_{i = 1}^n {A_i^2{B_i}} \sum\limits_{i = 1}^n {{A_i}{C_i}}  \le \sum\limits_{i = 1}^n {A_i^2} \sum\limits_{i = 1}^n {{A_i}{B_i}{C_i}} .\]

\textbf{Lemma 2.7}(see [1]) Let ${\lambda _i}$, $i = 1,2, \cdot
\cdot \cdot ,$ be the eigenvalues of (1.1), and ${u_i}$ the
corresponding eigenfunctions, then
\[\int_\Omega  {{u_i}{L^{k}}{u_i}}=\int_\Omega  {{{\left| {{\nabla ^k}{u_i}} \right|}^2}}  \le \left(\int_\Omega  {{u_i}{L^{l}}{u_i}}\right)^{\frac{k}{l}}=\lambda _i^{\frac{k}{l}}, k = 1, \cdot  \cdot  \cdot ,l - 1, \]
where
\[{\nabla ^k} \equiv \left\{ \begin{array}{l}
{\Delta ^{\frac{k}{2}}},{\kern 1pt} {\kern 1pt} {\kern 1pt} {\kern
1pt} {\kern 1pt} {\kern 1pt} {\kern 1pt} {\kern 1pt} {\kern 1pt}
{\kern 1pt} {\kern 1pt} {\kern 1pt} {\kern 1pt} {\kern 1pt} {\kern
1pt} {\kern 1pt} {\kern 1pt} {\kern 1pt} {\kern 1pt} {\kern 1pt}
{\kern 1pt} {\kern 1pt} {\kern 1pt} {\kern 1pt} {\kern 1pt} {\kern
1pt} {\kern 1pt} {\kern 1pt} {\kern 1pt} {\kern 1pt} {\kern 1pt}
{\kern 1pt} {\kern 1pt} {\kern 1pt} {\kern 1pt} {\kern 1pt} {\kern
1pt} {\kern 1pt} {\kern 1pt} {\rm{if  }}{\kern 1pt} ~k{\rm{ ~is
~even}},
\\
\nabla \left( {{\Delta ^{\frac{{k - 1}}{2}}}} \right),~{\rm{if
}}~k{\rm{ ~is ~odd}}.
\end{array} \right.\]


\section{Proofs of results}

\label{3} \textbf{Proof of Theorem 1.1. } We apply (2.3) with $A =
{L^l} = {\left( { - \Delta } \right)^l},{B_1} = {x_1}, \cdot  \cdot
\cdot ,{B_n} = {x_n},{T_1} = \frac{\partial }{{\partial {x_1}}},
\cdot  \cdot  \cdot ,{T_n} = \frac{\partial }{{\partial {x_n}}},f(x)
= {(\lambda  - x)^\alpha },g(x) = {(\lambda  - x)^\beta },$ and
obtain
\begin{align}
&{\left( {\sum\limits_{i = 1}^k {\sum\limits_{p = 1}^n {{{({\lambda
_{k + 1}} - {\lambda _i})}^\alpha }{{\left\langle {\left[
{\frac{\partial }{{\partial {x_p}}},{x_p}} \right]{u_i},{u_i}}
\right\rangle }_{{L^2}}}} } }
\right)^2}\notag
\\&\le 4\left( {\sum\limits_{i = 1}^k {\sum\limits_{p = 1}^n {{{({\lambda _{k + 1}} - {\lambda _i})}^\beta }{{\left\langle {\left[ {{L^l},{x_p}} \right]{u_i},{x_p}{u_i}} \right\rangle }_{{L^2}}}} } }
\right)\notag
\\&\times \left( {\sum\limits_{i = 1}^k {\sum\limits_{p = 1}^n {{{({\lambda _{k + 1}} - {\lambda _i})}^{2\alpha  - \beta  - 1}}\left\| {\frac{\partial }{{\partial {x_p}}}{u_i}} \right\|_{{L^2}}^2} } }
\right).
\end{align}
Since
\[\left[ {\frac{\partial }{{\partial {x_p}}},{x_p}} \right]{u_i} = \frac{\partial }{{\partial {x_p}}}({x_p}{u_i}) - {x_p}\frac{\partial }{{\partial {x_p}}}{u_i}= {u_i},\]
and
\[{\left\langle {\left[ {\frac{\partial }{{\partial {x_p}}},{x_p}} \right]{u_i},{u_i}} \right\rangle _{{L^2}}} = 1,\]
it arrives at \begin{equation}
{\left( {\sum\limits_{i = 1}^k
{\sum\limits_{p = 1}^n {{{({\lambda _{k + 1}} - {\lambda
_i})}^\alpha }{{\left\langle {\left[ {\frac{\partial }{{\partial
{x_p}}},{x_p}} \right]{u_i},{u_i}} \right\rangle }_{{L^2}}}} } }
\right)^2} = {\left( {n\sum\limits_{i = 1}^k {{{({\lambda _{k + 1}}
- {\lambda _i})}^\alpha }} } \right)^2}.
\end{equation}
By (2.1) and (2.2), it follows
\[\begin{array}{l}
{\left\langle {\left[ {{L^l},{x_p}} \right]{u_i},{x_p}{u_i}} \right\rangle _{{L^2}}}\\
 =  - 2l\int_\Omega  {{x_p}{u_i}{L^{l - 1}}\frac{\partial }{{\partial {x_p}}}{u_i}}  =  - 2l\int_\Omega  {\frac{\partial }{{\partial {x_p}}}{u_i}{L^{l - 1}}({x_p}{u_i})} \\
 =  - 2l\int_\Omega  {\frac{\partial }{{\partial {x_p}}}{u_i}\left\{ {{x_p}{L^{l - 1}}{u_i} - 2(l - 1){L^{l - 2}}\frac{\partial }{{\partial {x_p}}}{u_i}} \right\}} \\
 = 2l\int_\Omega  {{x_p}{u_i}{L^{l - 1}}\frac{\partial }{{\partial {x_p}}}{u_i}}  + 2l\int_\Omega  {{u_i}{L^{l - 1}}{u_i}}  - 4l(l - 1)\int_\Omega  {{u_i}{L^{l - 2}}\frac{{{\partial ^2}}}{{\partial
 x_p^2}}{u_i}},
\end{array}\]
hence
\[{\left\langle {\left[ {{L^l},{x_p}} \right]{u_i},{x_p}{u_i}} \right\rangle _{{L^2}}} = l\int_\Omega  {{u_i}{L^{l - 1}}{u_i}}  - 2l(l - 1)\int_\Omega  {{u_i}{L^{l - 2}}\frac{{{\partial ^2}}}{{\partial x_p^2}}{u_i}} .\]
We see from Lemma 2.7 that
\[\begin{array}{l}
\sum\limits_{p = 1}^n {{{\left\langle {\left[ {{L^l},{x_p}} \right]{u_i},{x_p}{u_i}} \right\rangle }_{{L^2}}}} \\
 = l(2l + n - 2)\int_\Omega  {{u_i}{L^{l - 1}}{u_i}} \\
 \le l(2l + n - 2){\left( {\int_\Omega  {{u_i}{L^l}{u_i}} } \right)^{\frac{{l - 1}}{l}}}\\
 = l(2l + n - 2)\lambda _i^{\frac{{l - 1}}{l}}
\end{array}\]
and \begin{equation}\sum\limits_{i = 1}^k {\sum\limits_{p = 1}^n
{{{({\lambda _{k + 1}} - {\lambda _i})}^\beta }{{\left\langle
{\left[ {{L^l},{x_p}} \right]{u_i},{x_p}{u_i}} \right\rangle
}_{{L^2}}}} } \le l(2l + n - 2)\sum\limits_{i = 1}^k {{{({\lambda
_{k + 1}} - {\lambda _i})}^\beta }\lambda _i^{\frac{{l - 1}}{l}}}.
\end{equation}
Since
\[\sum\limits_{p = 1}^n {\left\| {\frac{\partial }{{\partial {x_p}}}{u_i}} \right\|_{{L^2}}^2}  = \int_\Omega  {{u_i}L{u_i}}  \le {\left( {\int_\Omega  {{u_i}{L^l}{u_i}} } \right)^{\frac{1}{l}}} = \lambda _i^{\frac{1}{l}},\]
it yields \begin{equation}\sum\limits_{i = 1}^k {\sum\limits_{p =
1}^n {{{({\lambda _{k + 1}} - {\lambda _i})}^{2\alpha  - \beta  -
1}}\left\| {\frac{\partial }{{\partial {x_p}}}{u_i}}
\right\|_{{L^2}}^2} }  \le \sum\limits_{i = 1}^k {{{({\lambda _{k +
1}} - {\lambda _i})}^{2\alpha  - \beta  - 1}}\lambda
_i^{\frac{1}{l}}}.
\end{equation}
Instituting (3.2), (3.3) and (3.4) into (3.1), we have
\[{\left( {n\sum\limits_{i = 1}^k {{{({\lambda _{k + 1}} - {\lambda _i})}^\alpha }} } \right)^2} \le 4l(2l + n - 2)\left( {\sum\limits_{i = 1}^k {{{({\lambda _{k + 1}} - {\lambda _i})}^\beta }\lambda _i^{\frac{{l - 1}}{l}}} } \right) \left( {\sum\limits_{i = 1}^k {{{({\lambda _{k + 1}} - {\lambda _i})}^{2\alpha  - \beta  - 1}}\lambda _i^{\frac{1}{l}}} } \right)\]
and Theorem 1.1 is proved.

\textbf{Proof of Corollary 1.1. } To obtain (1.3), it suffices to
take $\alpha  = \beta  = 2$ in (1.2).

\textbf{Proof of Corollary 1.2. } When $1 \le \alpha  = \beta  \le
2$, we have from (1.2) that
\[{\left( {n\sum\limits_{i = 1}^k {{{({\lambda _{k + 1}} - {\lambda _i})}^\alpha }} } \right)^2} \le 4l(2l + n - 2)\left( {\sum\limits_{i = 1}^k {{{({\lambda _{k + 1}} - {\lambda _i})}^\beta }\lambda _i^{\frac{{l - 1}}{l}}} } \right) \left( {\sum\limits_{i = 1}^k {{{({\lambda _{k + 1}} - {\lambda _i})}^{2\alpha  - \beta  - 1}}\lambda _i^{\frac{1}{l}}} } \right).\]
Applying Lemma 2.5 to $\left( {\sum\limits_{i = 1}^k {{{({\lambda
_{k + 1}} - {\lambda _i})}^\beta }\lambda _i^{\frac{{l - 1}}{l}}} }
\right)$ and $\left( {\sum\limits_{i = 1}^k {{{({\lambda _{k + 1}} -
{\lambda _i})}^{2\alpha  - \beta  - 1}}\lambda _i^{\frac{1}{l}}} }
\right)$, it follows
\begin{align}
&{\left( {\sum\limits_{i = 1}^k {{{({\lambda _{k + 1}} - {\lambda
_i})}^\alpha }} } \right)^2}\notag
\\& \le \frac{{4l(2l + n - 2)}}{{{n^2}{k^2}}}\left( {\sum\limits_{i = 1}^k {{{({\lambda _{k + 1}} - {\lambda _i})}^\beta }} } \right)  \left( {\sum\limits_{i = 1}^k {{{({\lambda _{k + 1}} - {\lambda _i})}^{2\alpha  - \beta  - 1}}} } \right)\left( {\sum\limits_{i = 1}^k {\lambda _i^{\frac{{l - 1}}{l}}} } \right)\left( {\sum\limits_{i = 1}^k {\lambda _i^{\frac{1}{l}}} }
\right)\notag\\
& = \frac{{4l(2l + n - 2)}}{{{n^2}{k^2}}}\left( {\sum\limits_{i =
1}^k {{{({\lambda _{k + 1}} - {\lambda _i})}^\alpha }} } \right)
\left( {\sum\limits_{i = 1}^k {{{({\lambda _{k + 1}} - {\lambda
_i})}^{\alpha  - 1}}} } \right)\left( {\sum\limits_{i = 1}^k
{\lambda _i^{\frac{{l - 1}}{l}}} } \right)\left( {\sum\limits_{i =
1}^k {\lambda _i^{\frac{1}{l}}} } \right)\notag
\end{align}
and then
\[\sum\limits_{i = 1}^k {{{({\lambda _{k + 1}} - {\lambda _i})}^\alpha }}  \le \frac{{4l(2l + n - 2)}}{{{n^2}{k^2}}}\left( {\sum\limits_{i = 1}^k {{{({\lambda _{k + 1}} - {\lambda _i})}^{\alpha  - 1}}} } \right)\left( {\sum\limits_{i = 1}^k {\lambda _i^{\frac{{l - 1}}{l}}} } \right)\left( {\sum\limits_{i = 1}^k {\lambda _i^{\frac{1}{l}}} } \right),\]
so
\[\sum\limits_{i = 1}^k {{{({\lambda _{k + 1}} - {\lambda _i})}^{\alpha  - 1}}\left( {\left( {{\lambda _{k + 1}} - {\lambda _k}} \right) - \frac{{4l(2l + n - 2)}}{{{n^2}{k^2}}}\left( {\sum\limits_{i = 1}^k {\lambda _i^{\frac{{l - 1}}{l}}} } \right)\left( {\sum\limits_{i = 1}^k {\lambda _i^{\frac{1}{l}}} } \right)} \right)}  \le 0.\]
Since ${\lambda _i} \le {\lambda _k}$ for all $i \le k$, we have
(1.4).

\textbf{Proof of Corollary 1.3. } We have from Theorem 1.1 that
\[{\left( {\sum\limits_{i = 1}^k {{{({\lambda _{k + 1}} - {\lambda _i})}^\alpha }} } \right)^2} \le \frac{{4l(2l + n - 2)}}{{{n^2}}}\left( {\sum\limits_{i = 1}^k {{{({\lambda _{k + 1}} - {\lambda _i})}^\beta }\lambda _i^{\frac{{l - 1}}{l}}} } \right) \left( {\sum\limits_{i = 1}^k {{{({\lambda _{k + 1}} - {\lambda _i})}^{2\alpha  - \beta  - 1}}\lambda _i^{\frac{1}{l}}} } \right)\]
and show (1.5).

\textbf{Proof of Corollary 1.4. } Let us take $\alpha  = \beta  = 2$
in (1.5) to obtain (1.6).

\textbf{Proof of Corollary 1.5. } When $\alpha  = \beta  = 2$, we
know from Corollary 1.3 that \begin{equation} \sum\limits_{i = 1}^k
{({\lambda _{k + 1}} - {\lambda _i}} {)^\alpha } \le \frac{{4l(n +
2l - 2)}}{{{n^2}}}\left( {\sum\limits_{i = 1}^k {{{({\lambda _{k +
1}} - {\lambda _i})}^{\alpha  - 1}}{\lambda _i}} }
\right).\end{equation}
Using Lemmas 2.4 and 2.5, it implies
\[\sum\limits_{i = 1}^k {({\lambda _{k + 1}} - {\lambda _i}} {)^\alpha } \ge \frac{1}{{{k^{\alpha  - 1}}}}{\left( {\sum\limits_{i = 1}^k {({\lambda _{k + 1}} - {\lambda _i}} )} \right)^\alpha } \ge {\left( {\sum\limits_{i = 1}^k {({\lambda _{k + 1}} - {\lambda _i}} )} \right)^{\alpha  - 1}}\left( {{\lambda _{k + 1}} - {\lambda _k}} \right)\]
and \[\frac{{4l(n + 2l - 2)}}{{{n^2}}}\sum\limits_{i = 1}^k
{({\lambda _{k + 1}} - {\lambda _i}} {)^{\alpha  - 1}}{\lambda _i}
\le \frac{{4l(n + 2l - 2)}}{{{n^2}k}}\left( {\sum\limits_{i = 1}^k
{({\lambda _{k + 1}} - {\lambda _i}} {)^{\alpha  - 1}}}
\right)\left( {\sum\limits_{i = 1}^k {{\lambda _i}} } \right),\]
then from (3.5) that
\[{\left( {\sum\limits_{i = 1}^k {({\lambda _{k + 1}} - {\lambda _i}} )} \right)^{\alpha  - 1}}\left( {{\lambda _{k + 1}} - {\lambda _k}} \right) \le \frac{{4l(n + 2l - 2)}}{{{n^2}k}}\left( {\sum\limits_{i = 1}^k {{{({\lambda _{k + 1}} - {\lambda _i})}^{\alpha  - 1}}} } \right)\left( {\sum\limits_{i = 1}^k {{\lambda _i}} } \right).\]
Since
\[{\left( {\sum\limits_{i = 1}^k {({\lambda _{k + 1}} - {\lambda _i}} )} \right)^{\alpha  - 1}} \ge \sum\limits_{i = 1}^k {{{({\lambda _{k + 1}} - {\lambda _i})}^{\alpha  - 1}}},\]
we have (1.7).

\textbf{Proof of Corollary 1.6. } When $1 \le \alpha  = \beta  \le
2$, we have from Corollary 1.3 that
\[\sum\limits_{i = 1}^k {({\lambda _{k + 1}} - {\lambda _i}} {)^\alpha } \le \frac{{4l(n + 2l - 2)}}{{{n^2}}}\sum\limits_{i = 1}^k {{{({\lambda _{k + 1}} - {\lambda _i})}^{\alpha  - 1}}{\lambda _i}} ,\]
then
\[{\lambda _{k + 1}}\sum\limits_{i = 1}^k {{{({\lambda _{k + 1}} - {\lambda _i})}^{\alpha  - 1}}}  - \sum\limits_{i = 1}^k {{{({\lambda _{k + 1}} - {\lambda _i})}^{\alpha  - 1}}{\lambda _i}}  \le \frac{{4l(n + 2l - 2)}}{{{n^2}}}\sum\limits_{i = 1}^k {({\lambda _{k + 1}} - {\lambda _i}} {)^{\alpha  - 1}}{\lambda _i},\]
i.e.,
\[{\lambda _{k + 1}}\sum\limits_{i = 1}^k {{{({\lambda _{k + 1}} - {\lambda _i})}^{\alpha  - 1}}}  \le \left( {1 + \frac{{4l(n + 2l - 2)}}{{{n^2}}}} \right)\sum\limits_{i = 1}^k {({\lambda _{k + 1}} - {\lambda _i}} {)^{\alpha  - 1}}{\lambda _i}.\]
By Lemma 2.5, it follows
\begin{align}
&\left( {1 + \frac{{4l(n + 2l - 2)}}{{{n^2}}}} \right)\sum\limits_{i
= 1}^k {({\lambda _{k + 1}} - {\lambda _i}} {)^{\alpha  -
1}}{\lambda _i}\notag
\\
 &\le \left( {1 + \frac{{4l(n + 2l - 2)}}{{{n^2}}}} \right)\frac{1}{k}\left( {\sum\limits_{i = 1}^k {({\lambda _{k + 1}} - {\lambda _i}} {)^{\alpha  - 1}}} \right)\left( {\sum\limits_{i = 1}^k {{\lambda _i}} }
\right)\notag
\end{align}
and
\[\left( {{\lambda _{k + 1}} - \left( {1 + \frac{{4l(n + 2l - 2)}}{{{n^2}}}} \right)\frac{1}{k}\left( {\sum\limits_{i = 1}^k {{\lambda _i}} } \right)} \right)\left( {\sum\limits_{i = 1}^k {{{({\lambda _{k + 1}} - {\lambda _i})}^{\alpha  - 1}}} } \right) \le 0.\]
Since $\left( {\sum\limits_{i = 1}^k {{{({\lambda _{k + 1}} -
{\lambda _i})}^{\alpha  - 1}}} } \right) \ge 0$, it implies (1.8).


Acknowledgements. This work was supported by the National Natural
Science Foundation of China (Grant Nos. 11271299, 11001221) and
Natural Science Foundation Research Project of Shaanxi Province
(2012JM1014).


\end{document}